\documentclass{amsart}

\newtheorem{theorem}{Theorem}[section]

\newtheorem{lemma}[theorem]{Lemma}

\newtheorem{corollary}[theorem]{Corollary}

\newcommand{\RR}{{\mathbb R}}
\newcommand{\CC}{{\mathbb C}}
\newcommand{\TT}{{\mathbb T}}

\newcommand{\DD}{{\mathbb D}}

 \title[Remarks on a paper of Geronimo and Johnson]
{Remarks on a paper of Geronimo and Johnson}

\author[F. Peherstorfer and P. Yuditskii]
{ F. Peherstorfer$^1$ and P. Yuditskii$^2$}
\thanks{$^{1,2}$ The work was supported by the Austrian Science
Found FWF, project number: P16390--N04}

\begin{document}

\maketitle

\subsection{Character--automorphic Hardy Spaces}
Let $E$ be a finite union of (necessary non--degenerate) arcs
on the unit circle $\TT$. The domain $\overline\CC\setminus E$
is conformally equivalent to the quotient of the
unit disk by the action of a discrete group $\Gamma=\Gamma(E)$.
Let $z:\DD\to \overline\CC\setminus E$ be a covering map,
$z\circ\gamma=z$, $\forall \gamma\in\Gamma$. In what follows we
assume the following normalization to be hold
$$
z:(-1,1)\to (a_0,b_0)\subset \TT\setminus E,
$$
where $(a_0,b_0)$ is a fixed gap,
$\TT\setminus E=\cup_{j=0}^g(a_j,b_j)$. In this case one can chose
a fundamental domain $\mathfrak F$ and a system of generators
$\{\gamma_j\}_{j=1}^g$ of $\Gamma$ such that
they are symmetric with respect to the complex conjugation:
$$
\overline{\mathfrak F}=\mathfrak F, \quad
\overline{\gamma_j}=\gamma_j^{-1}.
$$

Denote by $\zeta_0\in\mathfrak F$ the preimage of the origin, 
$z(\zeta_0)=0$, then $z(\overline{\zeta_0})=\infty$.
Let $B(\zeta,\zeta_0)$ and $B(\zeta,\overline{\zeta_0})$
be the Green functions with
$B(\overline{\zeta_0},\zeta_0)>0$ and 
$B(\zeta_0,\overline{\zeta_0})>0$. Then
\begin{equation}\label{f1}
z(\zeta)=e^{ic}\frac{B(\zeta,{\zeta_0})}{B(\zeta,\overline{\zeta_0})}.
\end{equation}
It is convenient to rotate (if necessary) the set $E$ and to 
think that $c=0$. Note that $B(\zeta,{\zeta_0})$ is a
character--automorphic function
$$
B(\gamma(\zeta),{\zeta_0})=\mu(\gamma)
B(\zeta,{\zeta_0}), \quad\gamma\in\Gamma,
$$
with a certain $\mu\in\Gamma^*$. By \eqref{f1}
$$
B(\gamma(\zeta),\overline{\zeta_0})=\mu(\gamma)
B(\zeta,\overline{\zeta_0}), \quad\gamma\in\Gamma.
$$

Recall that the space $A^2_1(\alpha)$, $\alpha\in\Gamma^*$ is formed
by functions of Smirnov class in $\DD$ such that
$$
f|[\gamma](\zeta):=
\frac{f(\gamma(\zeta))}{\gamma_{21}\zeta+\gamma_{22}}
=\alpha(\gamma)f(\zeta),
\quad
\gamma=\begin{bmatrix} \gamma_{11}&\gamma_{12}\\
\gamma_{21}&\gamma_{22}
\end{bmatrix},
$$
and
$$
||f||^2:=\int_{\TT/\Gamma}|f(t)|^2\,dm(t)<\infty.
$$
We denote by $k^{\alpha}(\zeta,\zeta_0)$ the reproducing kernel
of this space and put
$$
K^{\alpha}(\zeta,\zeta_0):=
\frac{k^{\alpha}(\zeta,\zeta_0)}{||k||}=
\frac{k^{\alpha}(\zeta,\zeta_0)}{
\sqrt{k^{\alpha}(\zeta_0,\zeta_0)}}.
$$
Notice that in our case
$f(\zeta)\in A^2_1(\alpha)$ implies
$\overline{f(\bar\zeta)}\in A^2_1(\alpha)$ and therefore
$$
K^{\alpha}(\zeta_0,\zeta_0)=
K^{\alpha}(\overline{\zeta_0},\overline{\zeta_0}).
$$

\subsection{A recurrence relation for reproducing kernels}
We start with 
\begin{theorem}
Systems
$$
\{K^{\alpha}(\zeta,\zeta_0),
B(\zeta,\zeta_0)K^{\alpha\mu^{-1}}(\zeta,\overline{\zeta_0})\}
$$
and
$$
\{K^{\alpha}(\zeta,\overline{\zeta_0}),
B(\zeta,\overline{\zeta_0})K^{\alpha\mu^{-1}}(\zeta,{\zeta_0})\}
$$
form orthonormal bases in the two dimensional space
spanned by $K^{\alpha}(\zeta,\zeta_0)$
and
$K^{\alpha}(\zeta,\overline{\zeta_0})$.
Moreover
\begin{equation}\label{f2}
\begin{matrix}
K^{\alpha}(\zeta,\overline{\zeta_0})=&
a(\alpha)K^{\alpha}(\zeta,{\zeta_0})
+\rho(\alpha)
B(\zeta,\zeta_0)K^{\alpha\mu^{-1}}(\zeta,\overline{\zeta_0}),
\\
K^{\alpha}(\zeta,{\zeta_0})=&\overline{a(\alpha)}
K^{\alpha}(\zeta,\overline{\zeta_0})
+\rho(\alpha)
B(\zeta,\overline{\zeta_0})K^{\alpha\mu^{-1}}(\zeta,{\zeta_0}),
\end{matrix}
\end{equation}
where
$$
a(\alpha)=a=\frac{K^{\alpha}(\zeta_0,\overline{\zeta_0})}
{K^{\alpha}(\zeta,{\zeta_0})},\quad
\rho(\alpha)=\rho=\sqrt{1-|a|^2}.
$$
\end{theorem}
\begin{proof}
Let us prove the first relation in \eqref{f2}.
It is evident that the vectors
$K^{\alpha}(\zeta,\zeta_0)$ and
$B(\zeta,\zeta_0)K^{\alpha\mu^{-1}}(\zeta,\overline{\zeta_0})$
are orthogonal, normalized and orthogonal
to all functions $f$ from $A^2_1(\alpha)$ such that
$f(\zeta_0)=f(\overline{\zeta_0})=0$, that is to functions that
form orthogonal compliment to the vectors
$K^{\alpha}(\zeta,\zeta_0)$
and
$K^{\alpha}(\zeta,\overline{\zeta_0})$. Thus
$$
K^{\alpha}(\zeta,\overline{\zeta_0})=
c_1 K^{\alpha}(\zeta,{\zeta_0})
+c_2
B(\zeta,\zeta_0)K^{\alpha\mu^{-1}}(\zeta,\overline{\zeta_0}).
$$
Putting $\zeta=\zeta_0$ we get $c_1=a$. Due to orthogonality
we have
$$
1=|a|^2+|c_2|^2.
$$
Now, put $\zeta=\overline{\zeta_0}$.
Taking into account that
$K^{\alpha}(\zeta_0,\overline{\zeta_0})
=\overline{K^{\alpha}(\overline{\zeta_0},\zeta_0)}$
and
$B(\overline{\zeta_0},\zeta_0)>0$ we prove that
$c_2$ being positive is equal to $\sqrt{1-|a|^2}$.

Note that simultaneously we proved that
$$
\rho=B(\overline{\zeta_0},\zeta_0)
\frac{K^{\alpha\mu^{-1}}(\overline{\zeta_0},\overline{\zeta_0})}
{K^{\alpha}(\overline{\zeta_0},\overline{\zeta_0})}.
$$
\end{proof}

\begin{corollary} A recurrence relation for reproducing kernels
generated by
the shift of $\Gamma^*$ on the character $\mu^{-1}$ is of the form
\begin{equation}\label{rr}
\begin{split}
B({\zeta},\zeta_0)&
\begin{bmatrix}
K^{\alpha\mu^{-1}}(\zeta,{\zeta_0}),
&-K^{\alpha\mu^{-1}}(\zeta,\overline{\zeta_0})
\end{bmatrix}
\\
=&
\begin{bmatrix}
K^{\alpha}(\zeta,{\zeta_0}),
&-K^{\alpha}(\zeta,\overline{\zeta_0})
\end{bmatrix}
\frac 1\rho\begin{bmatrix}
1&a\\
\bar a&1
\end{bmatrix}
\begin{bmatrix}
z&0\\
0&1
\end{bmatrix}.
\end{split}
\end{equation}
\end{corollary}
\begin{proof}
We write
\begin{equation*}
\begin{split}
B({\zeta},\zeta_0)&
\begin{bmatrix}
K^{\alpha\mu^{-1}}(\zeta,{\zeta_0}),
&-K^{\alpha\mu^{-1}}(\zeta,\overline{\zeta_0})
\end{bmatrix}\\
=&
\begin{bmatrix}B({\zeta},\overline{\zeta_0})
K^{\alpha\mu^{-1}}(\zeta,{\zeta_0}),
&-B({\zeta},\zeta_0)K^{\alpha\mu^{-1}}(\zeta,\overline{\zeta_0})
\end{bmatrix}
\begin{bmatrix}
z&0\\
0&1
\end{bmatrix}.
\end{split}
\end{equation*}
Then, use \eqref{f2}.
\end{proof}
\begin{corollary}\label{c3}
Let
\begin{equation}\label{sa}
s^\alpha(z):=\frac{K^{\alpha}(\zeta,\overline{\zeta_0})}
{K^{\alpha}(\zeta,{\zeta_0})}
\end{equation}
Then the Schur parameters of the function
$\tau s^\alpha(z)$, $\tau\in \TT$, are
$$
\{\tau a(\alpha\mu^{-n})\}_{n=0}^\infty.
$$
\end{corollary}

\begin{proof}
Let us note that \eqref{rr} implies
$$
s^\alpha(z)=
\frac{a(\alpha)+z s^{\alpha\mu^{-1}}(z)}
{1+\overline{a(\alpha)}z s^{\alpha\mu^{-1}}(z)}.
$$
Then we iterate this relation. Also, multiplication
by $\tau\in \TT$ of a Schur class function evidently
leads to multiplication by $\tau$ of all Schur parameters.
\end{proof}

\noindent{\bf Remark}.
Let
\begin{equation}\label{Mat}
M(z;\alpha,\tau)=\frac{1+z\tau s^\alpha(z)}{1-z\tau s^\alpha(z)},
\end{equation}
$(\alpha,\tau)\in \Gamma^*\times\TT\simeq\TT^{g+1}$.
Then
$$
M(z;\alpha,\tau)=\int\frac{t+z}{t-z}d\sigma(t;\alpha,\tau)
$$
gives $g+1$ parametric family of probabilistic measures on the
unit circle.
Let us point out the normalization conditions for $M$:
$M(0)=1, M(\infty)=-1$.

\subsection{Example (one--arc case)} In this case
$\overline\CC\setminus E\simeq \DD$,
$\Gamma$ is trivial, and
\begin{equation}
z=\frac{B(\zeta,\zeta_0)}{
B(\zeta,\overline{\zeta_0)}}=
\frac{\frac{\zeta-\zeta_0}{1-\zeta\overline{\zeta_0}}
\overline{\left(\frac{\overline{\zeta_0}-\zeta_0}{
1-\overline{\zeta_0}^2}\right)}}
{\frac{\zeta-\overline{\zeta_0}}{1-\zeta\zeta_0}
\overline{\left(\frac{\overline{\zeta_0}-\zeta_0}{
1-{\zeta_0}^2}\right)}}=
-\frac{\zeta-\zeta_0}{\zeta-\overline{\zeta_0}}
\frac{1-\zeta\zeta_0}{1-\zeta\overline{\zeta_0}}
\frac{1-\overline{\zeta_0}^2}{1-{\zeta_0}^2}.
\end{equation}
That is
$$
b_0=z(1)=-\frac{1-{\zeta_0}}{1+{\zeta_0}}
\frac{1+\overline{\zeta_0}}{1-\overline{\zeta_0}},
$$
and $a_0=\overline{b_0}$.
We can put $\zeta_0=ir$, $0<r<1$. Then
$$
b_0=\left(\frac{2r}{1+r^2}+i
\frac{1-r^2}{1+r^2}\right)^2=e^{2i\theta},
$$
where
$$
\sin\theta=\frac{1-r^2}{1+r^2},\quad\theta\in(0,\pi/2).
$$
Further, for such $z$
$$
s(z)=
\frac{K(\zeta,\overline{\zeta_0})}{K(\zeta,{\zeta_0})}=
\frac{\frac{1}{1-\zeta{\zeta_0}}}
{\frac{1}{1-\zeta\overline{\zeta_0}}}
=\frac{1-\zeta\overline{\zeta_0}}
{1-\zeta{\zeta_0}}.
$$
Thus
$$
a=s(0)=\frac{1-|\zeta_0|^2}{1-\zeta_0^2}=\frac{1-r^2}{1+r^2}
=\sin\theta.
$$
The Schur parameters of the function $s_\tau(z)=\tau s(z)$ are
$$
s_\tau(z)\sim\{\tau\sin\theta,
\tau\sin\theta,\tau\sin\theta...\}.
$$

\subsection{Lemma on the reproducing kernel} 
Let us map (the unit circle of) $z$--plane onto (the upper
half--plane of) $\lambda$--plane in such a way that
$a_0\mapsto 1$, $z(0)\mapsto \infty$, $b_0\mapsto -1$.
In this way ($\zeta\mapsto z\mapsto\lambda$) we get the function
$\lambda=\lambda(\zeta)$ such that
\begin{equation}\label{f5}
z=\frac{B(0,\zeta_0)}
{B(0,\overline{\zeta_0})}
\frac{\lambda-\lambda_0}{\lambda-\overline{\lambda_0}},
\quad \lambda_0:=\lambda(\zeta_0).
\end{equation}

\begin{lemma} Let 
$k^\alpha(\zeta)=k^\alpha(\zeta,0)$ and
$B(\zeta)=B(\zeta,0)$ be subject to the normalization
$(\lambda B)(0)>0$. 
Denote by $\mu_0$ the character generated by $B$, i.e.,
$B\circ\gamma=\mu_0(\gamma)B$.
Then
\begin{equation}\label{f6}
k^\alpha(\zeta,\zeta_0)=
(\lambda B)(0)\frac
{\overline{k^\alpha(\zeta_0)}
\frac{k^{\alpha\mu_0}(\zeta)}{B(\zeta)k^{\alpha\mu_0}(0)}
-\overline{
\frac{k^{\alpha\mu_0}(\zeta_0)}{B(\zeta_0)k^{\alpha\mu_0}(0)}}
k^\alpha(\zeta)}
{\lambda-\overline{\lambda_0}}.
\end{equation}
\end{lemma}
\begin{proof}
We start with the evident orthogonal decomposition
$$
A^2_1(\alpha\mu_0)=\{k^{\alpha\mu_0}\}
\oplus B A^2_1(\alpha).
$$
We use this decomposition to obtain
$$
\lambda B f=(\lambda B)(0)f(0)
\frac{k^{\alpha\mu_0}(\zeta)}{k^{\alpha\mu_0}(0)}
+B\tilde f,\quad
\tilde f\in A^2_1(\alpha).
$$
Dividing by $B$ and using the orthogonality of the
summands, we get
\begin{equation}\label{f7}
P_+(\alpha)\lambda f=\tilde f=
\lambda f-
(\lambda B)(0)f(0)
\frac{k^{\alpha\mu_0}(\zeta)}{B(\zeta)k^{\alpha\mu_0}(0)},
\end{equation}
where $P_+(\alpha)$ is the orthoprojector onto
$A^2_1(\alpha)$.

Thus, on the one hand, for arbitrary $f\in A^2_1(\alpha)$
\begin{equation}\label{f8}
\langle(\lambda-\lambda_0)f,k^\alpha(\zeta,\zeta_0)\rangle
=\{P_+(\alpha)(\lambda-\lambda_0)f\}(\zeta_0).
\end{equation}
By virtue of \eqref{f7} we have
\begin{equation}\label{f9}
\begin{split}
\{P_+(\alpha)(\lambda-\lambda_0)f\}(\zeta_0)
&=
\lambda(\zeta_0) f(\zeta_0)-(B\lambda)(0) f(0)
\frac{k^{\alpha\mu_0}(\zeta_0)}{B(\zeta_0)k^{\alpha\mu_0}(0)}
\\
-\lambda_0f(\zeta_0)&=
-(B\lambda)(0) 
\frac{k^{\alpha\mu_0}(\zeta_0)}{B(\zeta_0)k^{\alpha\mu_0}(0)}
\langle f,k^\alpha\rangle.
\end{split}
\end{equation}
On the other hand, since the function
$\lambda$ is real on $\TT$,
\begin{equation}\label{f10}
\begin{split}
\langle(\lambda-\lambda_0)f,k^\alpha(\zeta,\zeta_0)\rangle
&=\langle
f,(\lambda-\overline{\lambda_0})k^\alpha(\zeta,\zeta_0)\rangle\\ 
&=\langle
f,P_+(\alpha)
(\lambda-\overline{\lambda_0})k^\alpha(\zeta,\zeta_0)\rangle.
\end{split}
\end{equation}
Comparing \eqref{f8} and \eqref{f9} with \eqref{f10},
we get
$$
P_+(\alpha)(\lambda-\overline{\lambda_0})k^\alpha(\zeta,\zeta_0)
=-\overline{(B\lambda)(0)
\frac{k^{\alpha\mu_0}(\zeta_0)}{B(\zeta_0)k^{\alpha\mu_0}(0)}
}k^\alpha(\zeta).
$$
Using \eqref{f7} again, we get
\begin{equation*}
\begin{split}
(\lambda-\overline{\lambda_0})k^\alpha(\zeta,\zeta_0)
-&
(B\lambda)(0) k^\alpha(0,\zeta_0)
\frac{k^{\alpha\mu_0}(\zeta)}{B(\zeta)k^{\alpha\mu_0}(0)}
\\
=&
-\overline{(B\lambda)(0)
\frac{k^{\alpha\mu_0}(\zeta_0)}{B(\zeta_0)k^{\alpha\mu_0}(0)}
}k^\alpha(\zeta).
\end{split}
\end{equation*}
Since $k^\alpha(0,\zeta_0)=\overline{k^\alpha(\zeta_0)}$, we have
\begin{equation*}
\begin{split}
&(\lambda-\overline{\lambda_0})k^\alpha(\zeta,\zeta_0)\\
&=(B\lambda)(0) 
\left\{\overline{k^\alpha(\zeta_0)}
\frac{k^{\alpha\mu_0}(\zeta)}{B(\zeta)k^{\alpha\mu_0}(0)}
-\overline{
\frac{k^{\alpha\mu_0}(\zeta_0)}{B(\zeta_0)k^{\alpha\mu_0}(0)}
}k^\alpha(\zeta)\right\}.
\end{split}
\end{equation*}
The lemma is proved.
\end{proof}

\begin{corollary}
In the introduced above notations
\begin{equation}\label{f11}
\begin{split}
z s^\alpha(z)
&=\frac{B(0,\zeta_0)}{B(0,\overline{\zeta_0})}
\frac{\lambda-\lambda_0}{\lambda-\overline{\lambda_0}}
\frac{K(\zeta,\overline{\zeta_0})}{K(\zeta,{\zeta_0})}\\
&=
\frac{B(0,\zeta_0)}{B(0,\overline{\zeta_0})}
\frac{{k^\alpha(\zeta_0)}
\frac{k^{\alpha\mu_0}(\zeta)}{B(\zeta)}
-{
\frac{k^{\alpha\mu_0}(\zeta_0)}{B(\zeta_0)}
}k^\alpha(\zeta)}{
\overline{k^\alpha(\zeta_0)}
\frac{k^{\alpha\mu_0}(\zeta)}{B(\zeta)}
-\overline{
\frac{k^{\alpha\mu_0}(\zeta_0)}{B(\zeta_0)}
}k^\alpha(\zeta)}
\\
&=
\frac{B(0,\zeta_0)k^\alpha(\zeta_0,0)}
{B(0,\overline{\zeta_0}){k^\alpha(0,\zeta_0)}}
\frac{r(\lambda;\alpha)-r(\lambda_0;\alpha)}
{r(\lambda;\alpha)-
\overline{r(\lambda_0;\alpha)}},
\end{split}
\end{equation}
where
\begin{equation}\label{ra}
r(\lambda;\alpha):=\frac{(\lambda B)(0)}{B(\zeta)}
\frac{k^{\alpha}(0)}{k^{\alpha\mu_0}(0)}
\frac{k^{\alpha\mu_0}(\zeta)}{k^{\alpha}(\zeta)}.
\end{equation}
\end{corollary}
Let us point out that 
functions \eqref{ra} are important in the 
spectral theory of Jacobi matrices [Sodin--Yuditskii],
they are
normalized by
$$
r(\lambda;\alpha)=\lambda+...,\quad\lambda\to \infty.
$$
\begin{corollary}
Let
$$
\tau(\alpha)=\left\{
\frac{B(0,\zeta_0)k^\alpha(\zeta_0,0)}
{B(0,\overline{\zeta_0}){k^\alpha(0,\zeta_0)}}\right\}^{-1}.
$$
Then
\begin{equation}\label{f12}
M(z;\alpha,\tau(\alpha))=
\frac{r(\lambda;\alpha)-\Re r(\lambda_0;\alpha)}
{i\Im r(\lambda_0;\alpha)}.
\end{equation}
\end{corollary}
\begin{proof}
By definition \eqref{Mat} and \eqref{f11}
$$
M(z;\alpha,\tau(\alpha))=
\frac{1+\frac{r(\lambda;\alpha)-r(\lambda_0;\alpha)}
{r(\lambda;\alpha)-\overline{r(\lambda_0;\alpha)}}}
{1-\frac{r(\lambda;\alpha)-r(\lambda_0;\alpha)}
{r(\lambda;\alpha)-\overline{r(\lambda_0;\alpha)}}}
=
\frac{r(\lambda;\alpha)-\Re r(\lambda_0;\alpha)}
{i\Im r(\lambda_0;\alpha)}.
$$
\end{proof}
\subsection{Main Theorem} Let $E$ be a finite union of arcs on the unit
circle $T$, $\TT\setminus E=\cup_{j=0}^g(a_j,b_j)$, normalized by the
condition $c=0$ in \eqref{f1}. Let $\mathfrak R(E)$ be the hyper--elliptic
Riemann surface with ramification points $\{a_j,b_j\}_{j=0}^g$
(double of the domain $\overline{\CC}\setminus E$).
Let us introduce a special collection of divisors on
$\mathfrak R(E)$:
\begin{equation*}
D(E)=\{D=\sum_{j=0}^g(t_j,\epsilon_j):
t_j\in [a_j,b_j],\quad \epsilon_j=\pm 1\},
\end{equation*}
where $(t_j,1)$ (correspondently  $(t_j,-1)$) denotes a point on the
upper (lower) sheet of the double $\mathfrak R(E)$, naturally,
$(a_j, 1)\equiv (a_j, -1)$ and $(b_j, 1)\equiv (b_j, -1)$.
Note that topologically $D(E)$ is the torus $\TT^{g+1}$.

Following [Akhiezer--Tomchuk, Pehersorfer--Steinbauer,
Geronimo--Johnson], we consider the
collection of functions 
$$
\mathfrak M(E)=\{M(z,D): D\in D(E)\}
$$
given in $\overline{\CC}\setminus E$ such that
$M(z,D)$ can be extended on $\mathfrak R(E)$ as a rational function
on it that has exactly $D$ as the divisor of poles and
meets the normalizations $M(0,D)=1$, $M(\infty,D)=-1$.
Note that the function is uniquely defined by $D$ and the normalizations
and has the integral representation
$$
M(z,D)=\int\frac{t+z}{t-z}\,d\sigma_D(t),
$$
with a probabilistic measure $\sigma_D$ on $\TT$.

\begin{theorem}
A given $D\in D(E)$ there exists a unique 
$(\alpha,\tau)\in \Gamma^*\times \TT$ such that the reflection
coefficients related to the orthogonal polynomials with respect
to $\sigma_D$ are $\{\tau a(\alpha\mu^{-n})\}_{n=0}^\infty$.
\end{theorem}
\begin{proof}
We only have to show that a given $M(z)=M(z,D)$ is of the form
\eqref{Mat} with  a certain  $(\alpha,\tau)$ and  to use 
Corollary \ref{c3}.

First, in the collection of functions
\begin{equation}\label{f16}
M_\theta(z)=
\frac{\cos\frac\theta 2 M(z) -i\sin\frac\theta 2}
{-i \sin\frac\theta 2 M(z)+\cos\frac\theta 2},
\quad 0\le\theta<2\pi,
\end{equation}
chose that one that has a pole at $z(0)\in (a_0,b_0)$.
Since $M(z)\in i\RR$ when $z\in (a_0,b_0)$, there exists a unique
$\theta$ that satisfied this condition. It is important, that
$M_\theta\in \mathfrak M(E)$, that is there exists a unique
$D_\theta$ such that
$M_\theta(z)=M(z,D_\theta)$.

Let us denote by $\tilde{\mathfrak R}(E)$ the Riemann surface 
that we obtain
by cutting and glueing two copies of the $\lambda$--plane 
(see \eqref{f5})
and by
$\tilde D$ the divisor on $\tilde{\mathfrak R}(E)$
that corresponds to a divisor $D\in D(E)$.
As it well known (see e.g. [Sodin-Yuditskii]), given $\tilde D_\theta$
there exists a unique $\alpha\in\Gamma^*$ 
such that $\tilde D_\theta$ is the divisor of poles
of a function of the form \eqref{ra}.

Now, consider the function
$$
M(z;\alpha,\tau(\alpha))=
\frac{r(\lambda(z);\alpha)-\Re r(\lambda_0;\alpha)}
{i\Im r(\lambda_0;\alpha)}
$$
with the chosen $\alpha$. It belongs to the class
$\mathfrak M(E)$ and, according to its definition,
has $D_\theta$ as the divisor of poles. Therefore, by uniqueness,
\begin{equation}\label{f17}
M_\theta(z)=M(z;\alpha,\tau(\alpha)).
\end{equation}
Substituting \eqref{f17} in \eqref{f16}
and solving for $M(z)$, we get
$$
M(z)=M(z;\alpha,\tau(\alpha) e^{-i\theta}).
$$
The theorem is proved.
\end{proof}

\end{document}